\documentclass{article}

\usepackage{amsmath}
\usepackage{amsfonts}
\usepackage{amssymb}
 \usepackage{amsthm}
\usepackage{graphicx}
\usepackage{array}
\usepackage{tabularx}
\usepackage{multirow}
\usepackage{indentfirst}
\usepackage{latexsym}
\usepackage{srcltx}
\usepackage{subfig}
\usepackage{wrapfig}

\newtheorem{theo}{Theorem}[section]

\newtheorem*{prof}{Proof}

\newtheorem{defi}[theo]{Definition}
\newtheorem{assumption}[theo]{Assumption}
\newtheorem{prop}[theo]{Proposition}

\newcommand{\NN}{\mathbb{N}}

\newcommand{\RR}{\mathbb{R}}

\newcommand{\Cc}{\mathcal{C}}
\newcommand{\Dd}{\mathcal{D}}

\newcommand{\Mm}{\mathcal{M}}

\newcommand{\Rr}{\mathcal{R}}

\renewcommand{\phi}{\varphi}

\newcommand{\Om}{\Omega}

\newcommand{\TTT}{T_{\alpha,\yd}}

\newcommand{\xtt}{\mathtt x}
\newcommand{\ytt}{\mathtt y}
\newcommand{\ztt}{\mathtt z}

\newcommand{\tu}{\tilde u}

\newcommand{\Cinf}{\text{C}^{\infty}}

\newcommand{\norm}[1]{|\!| #1 |\!|}

\renewcommand{\div}{\text{div}}

\newcommand{\abs}[1]{\left| #1 \right|}
\newcommand{\inner}[2]{\langle #1,#2 \rangle}

\DeclareMathOperator{\argmin}{argmin}

\DeclareMathOperator{\Tr}{Tr}

\newcommand{\yd}{y^\delta}
\newcommand{\uad}{u_\alpha^\delta}

\newcommand\set[1]{\left\{ #1 \right\}}

\title{Convergence of Variational Regularization Methods for Imaging on Riemannian Manifolds}

\author{{\large Nicolas Thorstensen${}^{1}$} and {\large Otmar Scherzer${}^{1,2}$}}

\date{}

\begin{document}

\maketitle

\begin{center}
{\normalsize
\begin{tabular}{cc}
${}^1$ Computational Science Center & ${}^2$ Radon Institute of Computational\\
 &  and Applied Mathematics\\
  University of Vienna & Austrian Academy of Sciences\\
  Nordbergstr. 15 & Altenberger Str.~69\\
  1090 Vienna, Austria & 4040 Linz, Austria\\
\end{tabular}}
\end{center}

\renewcommand{\thefootnote}{\fnsymbol{footnote}}
\renewcommand{\thefootnote}{\arabic{footnote}}

\begin{abstract}
We consider abstract operator equations $Fu=y$, where $F$ is a compact linear operator
between Hilbert spaces $U$ and $V$, which are function spaces on \emph{closed, finite dimensional Riemannian manifolds},
respectively. This setting is of interest in numerous applications such as Computer Vision and
non-destructive evaluation.

In this work, we study the approximation of the solution of the ill-posed operator equation with
Tikhonov type regularization methods. We prove well-posedness, stability, convergence, and convergence rates
of the regularization methods. Moreover, we study in detail the numerical analysis and
the numerical implementation. Finally, we provide for three different inverse problems numerical experiments.

\end{abstract}

\section{Introduction}
The problem of solving linear inverse and ill--posed problems has a long tradition in engineering (see \cite{Gro11}).
Several strategies have been proposed in the literature to solve such problems approximatively in a stable manner.

However, in most applications the data are assumed to be functions, which are defined on a subset of an Euclidean space.
In this paper the focus is on \emph{imaging} problems, where the data are functions on
\emph{closed, finite dimensional Riemannian manifolds}.
Such problems appear in Computer Vision and non-destructive evaluation, to name but a few (cf. Section \ref{sec:res}).

In this paper we take an abstract point of view and formulate the ill--posed imaging problem as the solution of
an operator equation
\begin{equation}
\label{eq:operator}
Fu=y\;.
\end{equation}
Here $F$ describes the physics of image formation, and $y$ denotes the ideal measurement data, which contains
neither noise $\delta$ nor modeling errors.
The operator $F: U\rightarrow V$ is a compact linear operator between Hilbert spaces of functions
defined on closed, finite dimensional Riemannian manifolds, respectively.
Consequently, the inverse operator is unbounded and the solution of (\ref{eq:operator}) is ill-posed.

In practice, ideal data are not available, but rather some approximation $\yd$. These perturbations, in general,
do not allow for a direct stable inversion of $F$. To provide a stable numerical solution, Tikhonov type regularization
is an adequate choice (see \cite{Gro84,Mor84,EngHanNeu96,TikLeoYag98,SchGraGroHalLen09}).
This method consisting in calculating a minimizer $\uad$ of the functional
\begin{equation}
\label{eq:solution}
\TTT(u) := \frac{1}{2} \norm{Fu-\yd}^{2} + \alpha \Rr(u)\,,
\end{equation}
which approximates the solution of (\ref{eq:operator}). Here, typically, $\Rr: U \to [0,+\infty]$ is a proper, convex regularization
functional. The parameter $\alpha$ controls the trade off between the quality of approximation of $Fu$ and $\yd$ as well as the
stability of the minimizer. The choice of the regularizing functional $\Rr$ is essential and is selected according to problem
specifications. Typical choices of $\Rr$, which are also considered in the paper, are Sobolev space (semi-)norms
and the total variation semi-norm on manifolds. In this paper, we do not consider more general settings of non-convex
regularization functionals, as it has been done in the Euclidean setting for instance in \cite{Zar09,Gra10}.

In this paper we consider three different applications of regularization methods on manifolds, which are \emph{denoising}, \emph{deblurring},
and an inverse problem from \emph{non-destructive evaluation}, which has been studied recently in \cite{LouRipSpiSpo11}.

In the following we summarize some related work: Diffusion filtering on surfaces has been used successfully for denoising
\cite{DiePreRumStr00,ClaDieRum04}, which can be considered a particular inverse problem (see Section \ref{sec:app}).
Even more multi-scale decomposition of data on manifolds can be used for denoising  \cite{AbrMouStaAfeBob07,GiaSloWen10}.
The numerical analysis and implementation of this paper is related to the work for discretization of partial differential equations on manifolds, in particular discretization of the Laplace-Beltrami operator on the manifold $\Mm$, $-\Delta_\Mm$. Here, in particular, we refer to pioneering work
of Dzuik \cite{Dzi88} on surface finite elements. The estimates there were generalized \cite{DemDzi07} by considering adaptive finite elements. Subsequently, several parabolic diffusion equations, \cite{EllDzi07a,ClaDieRum03} (isotropic) and \cite{ClaDieRum04} (anisotropic), were developed for manifold valued data. This topic should not be confused with the topic of the paper, where the domain of the functions is a manifold, where the context of the other papers is that the functions range in a manifold (see \cite{SapRin96}).

The paper is organized as follows:
In Section \ref{sec:analysis}, we prove well-posedness of variational regularization on closed, finite dimensional Riemannian manifolds. Also, convergence rates with respect to Bregman distances are obtained in the convex regularization setting, under a standard source condition. Section \ref{sec:res} is concerned with numerical minimization of the discrete Tikhonov functional - this is most probably the most important contribution of this paper. We provide a consistent discretization of convex Tikhonov functionals and formulate them in a purely matrix analysis fashion. As a byproduct this approach provides a consistent discretization of some nonlinear partial differential operators. Moreover, the consistent discretization is the basis to solve inverse problems in a stable way.
Section \ref{sec:app} provides numerical experiments for three different applications.
Finally, in Section \ref{sec:def}, we provide the basic notions on differential geometry and non-linear analysis on manifolds and provide some embedding results for Sobolev space and the space of functions of finite total variation.


\section{Analysis of variational regularization for functions on Riemannian manifolds}
\label{sec:analysis}
In this section, we state an analysis of variational regularization methods for solving the ill--posed operator Equation (\ref{eq:operator})
for functions on manifolds.
Well definedness, stability, convergence, and convergence rate are proven along the lines of \cite{SchGraGroHalLen09} - the manifold setting does
not further complicate the analysis, and thus is omitted. However, the results are formulated below for the sake of
completeness and fixation of the notation:
\begin{assumption}
\label{newassump}
\smallskip
\vspace{2mm}
\begin{enumerate}
\item[$(A1)$]  \label{it1} $U$ and $V$ are Hilbert spaces and $\tau_U$,$\tau_V$ denote the weak topologies, respectively.
\item[$(A2)$]  \label{it4} The functional $\Rr : U \to [0 ,+\infty]$ is convex and sequentially lower semi-continuous with respect to $\tau_U$.
\item [$(A3)$] \label{it5} $\Dd := \Dd (F) \cap \mathcal D ( \Rr) \neq \emptyset$ (which in particular implies that $\Rr$ is proper).
\item [$(A4)$] \label{it6} For every $\alpha > 0$ and $C>0$, the lower level set of the Tikhonov functional
\[
\mathrm{level}_{C}(\TTT):=\{u \in U:\TTT(u) \leq C \}
\]
are sequentially pre-compact with respect to $\tau_U$.
\item [$(A5)$]  \label{it7}  For every $\alpha > 0$ and $C>0$, the set $\mathrm{level}_{C}(\TTT)$ is sequentially closed with respect to $\tau_U$ and the restriction of $F$ to $\mathrm{level}_{C}(\TTT)$ is sequentially continuous with respect to $\tau_U$ and  $\tau_V$.
\end{enumerate}
\end{assumption}

The results from \cite{SchGraGroHalLen09} imply then:
\begin{theo}
\label{Well-posd}
Let Assumption \ref{newassump} hold. Then,
\begin{itemize}
\item there exists a minimizer of $\TTT$ for every $\alpha > 0$ and $\yd \in V$.
\item Let $\alpha > 0$. Then, for every sequence $y_k \rightarrow \yd$ let us denote
\[
u_k \in \argmin T_{\alpha,y_k}, \quad k \in \NN,
\]
then $(u_k)$ has a convergent subsequence. Every convergent subsequence of $(u_k)$ converges to a minimizer of $\TTT$.
\item Assume that Equation (\ref{eq:operator}) has a solution in $\Dd$. Moreover,
assume that a function $\alpha:(0,\infty) \rightarrow(0,\infty)$ satisfies
\[
\alpha(\delta) \rightarrow 0 \text{ and } \frac{\delta^2}{\alpha(\delta)} \rightarrow 0 \text{ for } \delta \rightarrow 0.
\]
Let the sequence $\delta_k$ of positive numbers converging to $0$ and assume that the data $y_k:= y^{\delta_{k}}$,
$\alpha_k := \alpha(\delta_k)$ satisfies
$\norm{y-y_k} \leq \delta_k$.

Then, $u_k \in \argmin T_{\alpha(\delta_k),v_k}$ has a convergent subsequence and every limit is a solution of
Equation (\ref{eq:operator}).
\end{itemize}
\end{theo}

For obtaining qualitative estimates for the convergence of a Tikhonov regularized solution to a minimum norm solution, some
additional assumptions, such as the so-called source condition, are needed.
\begin{prop}[Convergence rates]
\label{pr:rates}
Let Assumption \ref{newassump} hold. Assume there exists a $\Rr$ minimizing solution $u^\dagger \in \Dd(F) \cap \Dd(\Rr)$ of Equation (\ref{eq:operator}) and an element
\begin{equation}
\label{eq:sourceCond}
\xi \in \text{Ran}(F^*) \cap \partial {\Rr}(u^\dagger).
\end{equation}
Then, with the parameter choice $\alpha \sim  \delta$, we have
\[
D_\xi (\uad,u^\dagger)= O(\delta) \text{ and } \norm{F\uad-\yd} =O(\delta)\,,
\]
where $D_\xi$ denotes the Bregman distance, which is defined as follows:
\begin{equation*}
D_\xi(\tu,u) := \Rr(\tu)-\Rr(u)- \inner{\xi}{\tu-u}\,,\quad
u,\tu \in U\;.
\end{equation*}
Here $ \inner{\xi}{\tu-u}$ denotes the inner product on the Hilbert space $U$.
\end{prop}

Proposition \ref{pr:rates} applies for instance to total variation minimization
\[
T_{\alpha,\yd}(u) := \frac{1}{2} \int_\Mm (F u - \yd)^2 + \alpha \abs{D_\Mm u}(\Mm)\,,
\]
where
\begin{equation}
\label{eq:tv}
\begin{aligned}
& \abs{D_\Mm u}(\Mm) \\
:= &\sup \set{ \int_{\Mm} u  \div_\Mm X \ d\nu(g) : X \in C^\infty(\Mm,\RR^n) \text{ and }
\norm{X}_{\text{\upshape L}^\infty(\Mm,\RR^n)} \leq 1}.
\end{aligned}
\end{equation}
denoted the total variation of $u$ on the manifold and $\nabla_\Mm$ denotes the covariant derivative.
We choose the space $U,V = \text{\upshape L}^2(\Mm)$. Moreover, we assume that $F$ is continuous on $\text{\upshape L}^2(\Mm)$
with $\Dd(F) \cap \Dd(\Rr) \neq \emptyset$. $\Rr (u)= \abs{D_\Mm u}(\Mm)$ is the total variation semi-norm.
The verification of Assumption \ref{newassump} is similar to the Euclidean setting, and thus omitted.
However, the verification requires the Meyer-Serrin Theorem \ref{meyerSerr} and the Compactness Theorem \ref{th:compactness} for functions of Bounded Variation  $BV(\Mm)$ on manifolds. Using both theorems allows to shows that the
Poincar\'e inequality holds and from this follows that $\TTT$ is coercive (A.4). The Compactness Theorem is applied to verify (A.5).

Interpretation of the source condition and the convergence rates have been given in \cite{SchGraGroHalLen09} for function
defined on subsets of $\RR^n$, but are valid in the manifold setting in a completely analogous manner.

\section{Numerical results}
\label{sec:res}
In this section we discuss the implementation of variational regularization method for functions defined on manifolds.
Afterwards three inverse problems and numerical experiments are considered. The three applications are \emph{denoising},
\emph{deblurring}, and an \emph{inverse problem} for the Funk--Radon transform.

Now we discuss the numerical minimization of the discretized Tikhonov functional.

We assume that the closed Riemannian manifold $\Mm$ is approximated by a polyhedron $M$ represented
as $M =(\mathcal V,\mathcal T)$, with vertices $\mathcal V =\set{v_1,\dots,v_{K}} \in \RR^3$  and triangles
$\mathcal T =\set{T_1,\dots,T_L} \subset \mathcal V \times \mathcal V \times \mathcal V$.
The three components of a vertex $v_k$ are denoted by $v^j_k \text{ with } j=1,\dots,3$. Each triangle $T_i \in \mathcal T$ is defined by the set of indices $\mathtt t_i =\set{i_1,i_2,i_3 }_i$ of the
vertices $\set{ v_{i_1},v_{i_2},v_{i_3}}$, which are assumed to be counter-clockwise oriented.
In this section we only deal with the manifold $M$, and assume that $M$ is a sufficiently good approximation to $\Mm$
which justifies an identification. Consequently, also the metric $g$ and the surface measure on $M$, $\nu(g)$,
are also identified.

The polyhedral surfaces used in the numerical experiments below have been taken from the database \cite{AIMSHAPE}.
Each surface is closed, of genus zero, and consists of approximately $25000$ vertices. For genus zero surfaces a natural
parametrization is the sphere. Following \cite{GotGuShe03} imaging testdata $\yd$ on the manifold $M$ is generated by mapping
a given function with planar domain $\Om \subset \RR^2$ onto $M$ by making use of the spherical parametrization.

\subsubsection*{Polyhedral representation}
Each triangle $T_i$ is parameterized with respect to its vertices $v_k$, $v_j$, and $v_l$ by using barycentric coordinates
\begin{equation*}
\xtt_{k,i} (\gamma) = v_{k} + \zeta_1 (v_{j}-v_{k}) + \zeta_2 (v_{l}-v_{k})\,,
\end{equation*}
where
\begin{equation*}
\Gamma = \set{ \gamma = (\zeta_1,\zeta_2): \zeta_1 \in [0,1] \text{ and } \zeta_2 \in [0,1-\zeta_1 ] }\;.
\end{equation*}

We approximate the minimizer of the Tikhonov functional from (\ref{eq:solution}) by the minimizer of $\TTT$ on the finite dimensional space
of piecewise linear functions on the polyhedron $M$:
For $k \in \set{1,\ldots,K}$ let $(\phi_k)$ be the function, which is continuous on $M$, linear on each triangle $T_i$, $i=1,\ldots,L$, and
satisfies $\phi_k(v_k)=1$ and $\phi_k(v_s)=0$ if $s \neq k$. On each triangle $T_i$ we have exactly three such functions (for every vertex).
The set of piecewise linear functions is the linear span of the functions $\phi_k$:
\begin{equation*}
PL(M) :=  \set{\sum_{k} u_k \phi_k}\;.
\end{equation*}
From the definition of $\phi_k$ and $\xtt_{k,i}$ it follows that
\begin{equation}
\label{id_for_xtt}
\xtt = \sum_k v_k \phi_k (\xtt) \text{ for every } \xtt \in M.
\end{equation}
and
\begin{equation}
\label{id_for_phik}
\phi_k(\xtt_{k,i}(\gamma))=1-\zeta_1-\zeta_2\text{ for every } \gamma \in \Gamma\;.
\end{equation}
Minimization of the  Tikhonov functional is performed for $u \in PL(M)$ and we assume that the data $\yd \in PL(M)$ too.
Thus the functions $u$ over which we minimize and the data $\yd$ can be expressed via there series expansion:
\begin{equation}
\label{eq:u_und_vektor}
u(\xtt) = \sum_{k} u_k \phi_k(\xtt) \text{ and } \yd(\xtt) = \sum_{k} \yd_k \phi_k(\xtt)\;.
\end{equation}
The vectors of coefficients are denoted in boldface by $\mathbf u:= (u_k)$, $\mathbf \yd := (\yd_k)$, respectively.
The Jacobian of the parametrization of the manifold $M$ is the matrix
\begin{equation}
\label{eq:J}
\mathbf J = \left( \begin{array}{rccl}
                   \mathbf J_1 & 0 & \ldots & 0 \\
                   0 & \mathbf J_2 & \ddots & \vdots \\
                   \vdots & \ddots & \ddots & \vdots\\
                   0 &  0 & \ldots & \mathbf J_L
                   \end{array} \right) \in \RR^{2L \times 3L}
\end{equation}
with $L$ blocks
\[
\mathbf J_i=
\left[
\begin{array}{ccc}
 \frac{\partial \phi_k}{\zeta_1} & \frac{\partial \phi_j}{\zeta_1}  & \frac{\partial \phi_l}{\zeta_1}    \\
\frac{\partial \phi_k}{\zeta_2} & \frac{\partial \phi_j}{\zeta_2}  & \frac{\partial \phi_l}{\zeta_2}
\end{array}
\right] =
\left[
\begin{array}{ccc}
-1 &1&0  \\
-1&0&1
\end{array}
\right].
\]
Each submatrix $\mathbf J_i$ is the Jacobian of the parameterizations in the triangle $T_i$.
All vertices in $M$ are put in a block diagonal matrix $\mathbf V \in \RR^{3L \times 3L}$ with  $L$ blocks
\[
\mathbf V_i=
\left[
\begin{array}{ccc}
 v_k^1 & v_j^1 & v_l^1\\
 v_k^2 & v_j^2 & v_l^2  \\
 v_k^3 & v_j^3 & v_l^3\\
\end{array}
\right] \;.
\]
Each submatrix $\mathbf V_i$ stores the vertices from triangle $T_i$ ordered accordingly to the basis functions $\phi_k$.
The metric tensor $\mathbf G \in \RR^{2L \times 2L}$ on $M$ is block diagonal matrix with the $L$ diagonal blocks
\begin{equation*}
\mathbf G_i := (\mathbf V_i\mathbf J_i^T)^T (\mathbf V_i \mathbf J_i^T  ) \in \RR^{2 \times 2}\;.
\end{equation*}
Again, $\mathbf G_i$ is the metric tensor in a given triangle $T_i$.
Let $A_i$ denote the area of the triangle $T_i$, then the volume of the metric tensor satisfies:
\[
\sqrt{\abs{g}} = \sqrt{|\text{det}(\mathbf G_i)|} =  2A_i\,,
\]
Therefore the surface measure $d\nu(g)$ can be expressed in barycentric coordinates and the relation reads as follows:
\begin{equation*}
d\nu(g) = 2 A_i d\gamma.
\end{equation*}
Let $\mathbf {\tilde{ V}} \in \RR^{3L \times K}$ a matrix which encodes the connectivity of the manifold $M$. That is
\begin{equation*}
\mathbf{ \tilde V}_{ki}=
\begin{cases}
 1& \text{ if } i \in \mathtt t_{k'} \text{ where } k=3(k'-1)+j,\, j=1,2,3
 \\
0 &\text{else.}
\end{cases}
\end{equation*}
$\mathbf {\tilde{V}_i} \in \RR^{3 \times K}$ is a linear mapping assigning each triangle $T_i$ the indices of the three vertices.
Accordingly, the covariant derivative of $u \in PL(M)$, $\nabla_M u$, on the triangle $T_i$ is a constant vector and is given by
\begin{equation}
\label{eq:nablaM}
\begin{aligned}
\mathbf Z_i :=  \mathbf V_i  \mathbf J^T_i \mathbf G_i^{-1} \mathbf J_i \mathbf{ \tilde V_i} \mathbf u   \in \RR^{3 \times1}\;.
\end{aligned}
\end{equation}
The matrix
\begin{equation*}
\mathbf Z = \left(
\begin{array}{rcccl}
\mathbf Z_1 & 0 & \ldots & 0 & 0\\
0 & \mathbf Z_2 & 0 & \ldots & 0\\
\vdots & 0 & \ddots & 0 &\vdots\\
0 & 0 & \ldots & 0 & \mathbf Z_L
\end{array}
\right) \in \RR^{3L \times L}
\end{equation*}
consists of the gradient vectors of $u$ on each triangle of $M$. The matrix $\mathbf Z^T \mathbf Z$ is a positive semi-definite diagonal matrix.
Thus $(\mathbf Z^T \mathbf Z)^{p/2}$ is the matrix consisting of the $p/2$ powers of diagonal entries.

\subsection{Discretization of the Tikhonov Functional}
In the following we consider minimization of the discrete Tikhonov functional with functions defined on $PL(M)$.
The goal is to express the fit-to-data term and the regularization functional in dependence of the vector $\mathbf u$.

In all our test cases we have that $F:U \to V$, where $U$ and $V$ are function spaces defined on the same closed, finite dimensional
Riemannian manifold $\Mm$.
We assume that $Fu$ can be approximated by a piecewise linear function on the polyhedron $M$ (note that here both $Fu$ and $\Mm$ are
approximated).

The linear operator $F$ may not necessarily map onto piecewise linear functions and thus the elements of the range are again
approximated by the discrete operator
\begin{equation}
F^d u(\xtt) := \sum_{k=1}^{K} \sum_{j=1}^K F_{kj} u_j  \phi_k(\xtt) \sim F u(\xtt)\;.
\end{equation}
In the following, for the sake of simplicity of notation, we identify the discrete operator $F^d$ with the matrix $\mathbf F$
of coefficients. Moreover, we assume that the discretization is fine enough that we can identify $F$ and $F^d$ on $PL(M)$.

We use the following approximations for the fit-to-data term and the regularization functional:
\begin{itemize}
\item
Let the matrix $\mathbf A \in \RR^{K \times 3L}$ be defined by the areas $A_i$ of the triangles:
\begin{equation*}
A_{ij}=
\begin{cases}
 \sqrt{2A_i} &\text{ if } 3(i-1) <  j <3i+1,\\
0 &\text{else.}
\end{cases}
\end{equation*}
Then,
\begin{equation}
\label{eq:discrete_fit}
\begin{aligned}
\frac{1}{2} \norm{Fu-\yd}_{\text{\upshape L}^2(M)}^2 &= \frac{1}{2}  \int_{M} (Fu(x)-\yd(x))^2 d\nu_g(x)\\
&= \frac{1}{2} \sum_{T_i}  \int_{T_i}(F u(\xtt)- \yd(\xtt))^2 d\nu_g(\xtt)\\
&= \frac{1}{6} \sum_{T_i} 2A_i \sum_{j \in \mathtt t_i} ((\mathbf F \mathbf u)_{j} - \mathbf y_j^\delta)^2\\
&= \frac{1}{6} \norm{\mathbf A( \mathbf{ \tilde V}  \mathbf F \mathbf u- \mathbf{ \tilde V}   \mathbf \yd)}^2
\end{aligned}
\end{equation}

\item In the applications presented below the regularization functional is either the total variation semi-norm or the quadratic Sobolev semi-norm of the gradient. We evaluate these functionals for $u \in PL(M)$ on the polyhedron:
    \begin{equation*}
    u \in PL(M) \to \int_M |\nabla_M u|^p d\nu(g)\,, \quad p=1,2\;.
    \end{equation*}

Let the diagonal matrix $\mathbf {\tilde A} \in \RR^{L \times L}$ be defined by the areas $A_i$ of the triangles:
\[
\mathbf {\tilde A_{ii}} = 2 A_i \; .
\]
From the above considerations we find that
\begin{equation}
\label{eq:TVDiscrete}
\begin{aligned}
\Rr(u) &=\frac{1}{p} \int_M \abs{\nabla_M u}^{p} d \nu(g)\\
        &=\frac{1}{p} \sum_{T_i} 2 A_i \abs{( \nabla_M u)_i}^{p}\\
	   &=\frac{1}{p} \Tr(   \mathbf {\tilde A}  (\mathbf Z^T  \mathbf Z )^\frac{p}{2})
\end{aligned}
\end{equation}
\end{itemize}
Because we have that
\begin{equation*}
\frac{\partial \mathbf Z}{\partial \mathbf u} = \mathbf V \mathbf J^T \mathbf G^{-1} \mathbf J \mathbf {\tilde V}\,,
\end{equation*}
it follows that
the derivative of the discrete functional $\TTT$ (with $\Mm$ replaced by $M$) at $\mathbf u$ is given by
\begin{equation*}
\begin{aligned}
\frac{ \partial \TTT}{\partial \mathbf u}(\mathbf u)
&= \frac{1}{6} \mathbf F^T \mathbf{ \tilde V^T}  \mathbf A^T \mathbf A ( \mathbf{ \tilde V}  \mathbf F \mathbf u- \mathbf{ \tilde V}  \mathbf \yd) \\
& \qquad + \alpha \mathbf{\tilde A}( \mathbf Z^T \mathbf Z )^{(p-2)/2}  \mathbf{ \tilde V^T} \mathbf J^T    \mathbf G^{{-1}^T}   \mathbf J \mathbf V^T
                   ( \mathbf V \mathbf J^T \mathbf G^{-1} \mathbf J \mathbf{ \tilde V} \mathbf u)\;.
\end{aligned}
\end{equation*}
The formal derivative of the $d \Rr(u)$ is a discrete approximation of the differential operator
\[
- \div_M \left(\abs{\nabla_M u}^{p-2} \nabla_M u \right)\;.
\]
In particular, for $p=2$ we obtain a consistent approximation of the Laplace-Beltrami operator.

The optimality condition $d \TTT(\mathbf u,\mathbf \rho) = 0$ for all $\mathbf \rho$, can be solved with
a Landweber fixed point iteration:
\[
\mathbf u^{(k+1)}=  \mathbf u^{(k)} -\kappa \nabla \TTT(\mathbf u^{(k)}) \quad k=0,1,2,\ldots.
\]
Here $\kappa$ denotes the step size and is chosen to satisfy a stability criterion \cite{Gro84}. The algorithm is usually
terminated if the difference of the update $\norm{\mathbf u^{(k+1)}- \mathbf u^{(k)}}_{\infty}$  is below a given threshold
for the first time.
\subsection{Applications}
\label{sec:app}
\subsubsection*{Denoising of data on manifolds}
We consider \emph{denoising} of image data on a closed finite dimensional Riemannian manifold. The usual assumption is that
the data $\yd$ can be decomposed into a ideal image $u^\dagger$ and additive white noise $n_\sigma$,
with mean $0$ and variance $\sigma$. That is
\[
\yd = u^\dagger + n_\sigma \text{ with }\norm{n_\sigma}_{\text{\upshape L}^2(\Mm)} \leq \sigma\;.
\]
This corresponds to Equation (\ref{eq:operator}) where the operator $F$ is the identity. Thus denoising can be viewed
as an inverse problem.
The goal of denoising is to remove the noise component $n_\sigma$ from $\yd$ but at the same time preserve the visual
appearance of the clean image $u^\dagger$.

\begin{figure}
  \centering
 \subfloat[$u^\dagger$ = Lenna on $M$]{\includegraphics[scale=0.4]{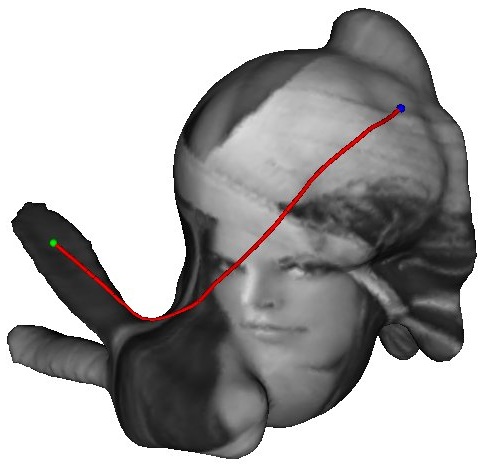}}
  \subfloat[Testdata: $u^\dagger + \mathcal N(0,\sigma)$ on $M$]{\includegraphics[scale=0.4]{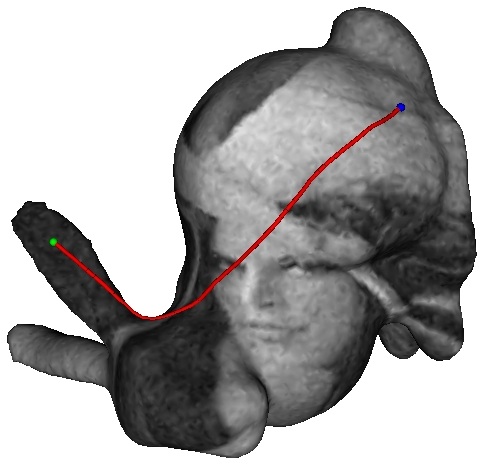}}

   \subfloat[Result with quadratic Sobolev semi-norm regularization]{\includegraphics[scale=0.4]{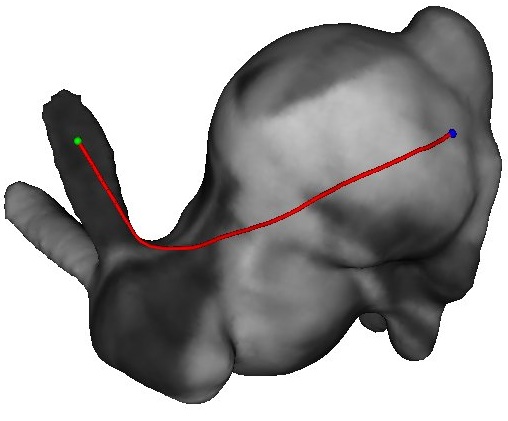}}
 \subfloat[$\uad-u^\dagger$]{\includegraphics[scale=0.4]{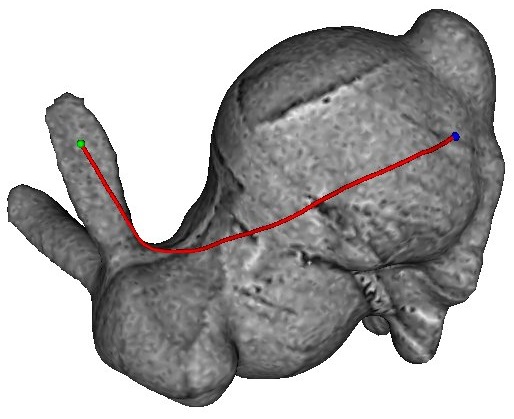}}

 \subfloat[Result with total variation denoising]{\includegraphics[scale=0.4]{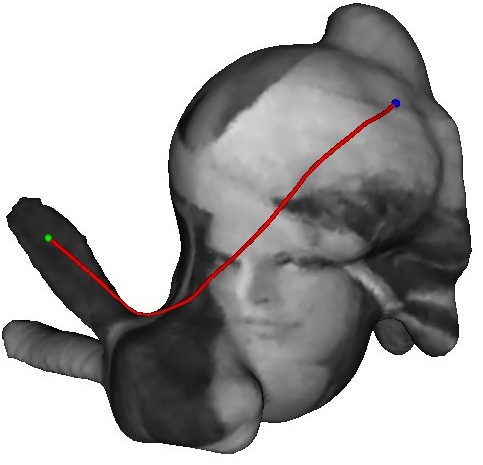}}
 \subfloat[$\uad-u^\dagger$]{\includegraphics[scale=0.4]{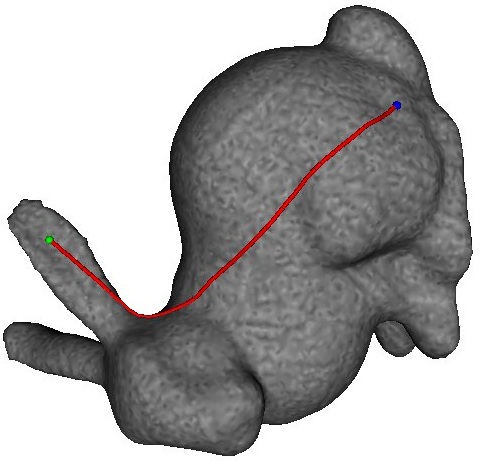}}
  \caption{\emph{Denoising:} The first row shows the ideal data $u$ and the noisy signal $\yd$. The second row depicts the results with total variation regularization and the residual image $\uad-u^\dagger$. The last rows results are obtained with by regularizing with the squared gradient.}
  \label{fig:resultsDenoising}
\end{figure}

In Figure \ref{fig:resultsDenoising} we compare quadratic Sobolev semi-norm regularization with total variation minimization.

\subsubsection*{Image deblurring}
The general assumption is that the imaging data $\yd$ is obtained from the clean image $u^\dagger$ by convolution with a smooth kernel
function $h$ and by additive white noise with mean zero and variance $\sigma$. Thus, in the terminology of the paper, the operator
equation reads as follows
\[
\yd = h \star u^\dagger + n_\sigma =: F u^\dagger + n_\sigma\;.
\]
In our numerical experiments we assume that the kernel function $h$ is a Gau\ss{}-function with variance $\tau$. That is,
\[
h(p,q) = \frac{1}{2 \pi \tau^2} \exp \left(\frac{d^2_g(p,q)}{2 \tau^2} \right) \text{ for every } p, q \in M\,,
\]
where $d_g(p,q)$ denotes the geodesic distance on the polyhedron $M$.

For implementing the Landweber algorithm (\ref{eq:steep_deconvolution}) we use the discrete convolution, which is
(similar as in Section \ref{sec:res}) written as a bold face matrix $\mathbf H$ with matrix entries
$d(v_{i,j},v_{k,l})$. The geodesic distance between two points $p$ and $q$ can be computed by solving the Eikonal equation with constant
velocity $\rho(\xtt)=1$ on $M$ with an algorithm described in \cite{KimSet98}.
That is, after fixing one point $q \in M$, $d(\xtt,q)$ solves the Eikonal equation:
\[
\abs{\nabla d(\xtt,q)} =  \rho(\xtt)\;.
\]

The Landweber algorithm for minimization of the discretized regularization functional reads as follows
\begin{equation}
\label{eq:steep_deconvolution}
\mathbf u^{(k+1)} = \mathbf u^{(k)} - \kappa (\mathbf H^T (\mathbf H \mathbf u^{(k)} - \mathbf \yd) + \alpha \nabla \Rr(\mathbf u^{(k)}))\;.
\end{equation}

\begin{figure}
  \centering
 \subfloat[$u^\dagger$ = Lenna on $M$]{\includegraphics[scale=0.4]{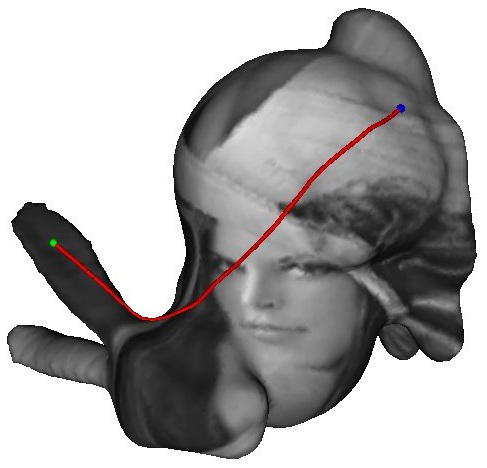}}
  \subfloat[$\yd$ = Gau\ss{}ian convolution of $u^\dagger$ + $\mathcal N(0,\sigma)$ ]{\includegraphics[scale=0.4]{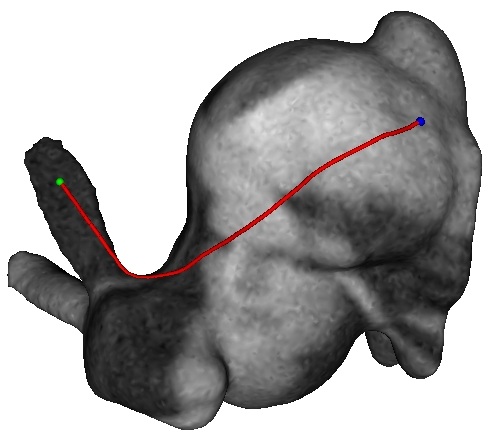}}

 \subfloat[Result with quadratic Sobolev semi-norm regularization]{\includegraphics[scale=0.4]{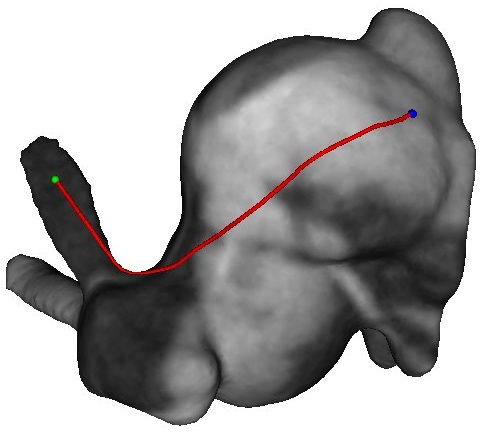}}
 \subfloat[$\uad-u^\dagger$]{\includegraphics[scale=0.4]{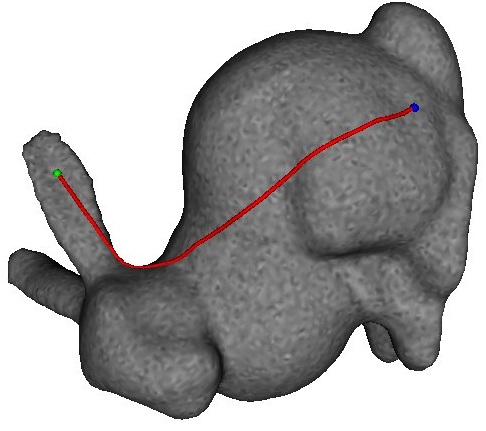}}

 \subfloat[Result with total variation regularization]{\includegraphics[scale=0.4]{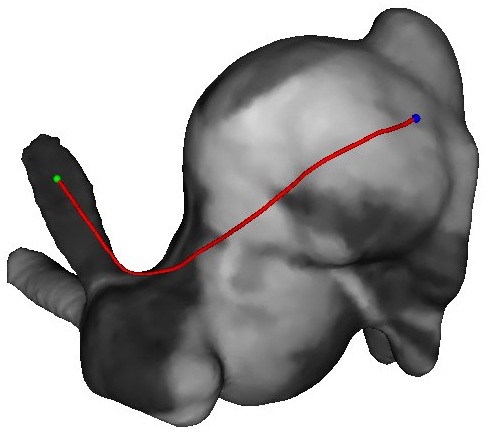}}
 \subfloat[$\uad-u^\dagger$]{\includegraphics[scale=0.4]{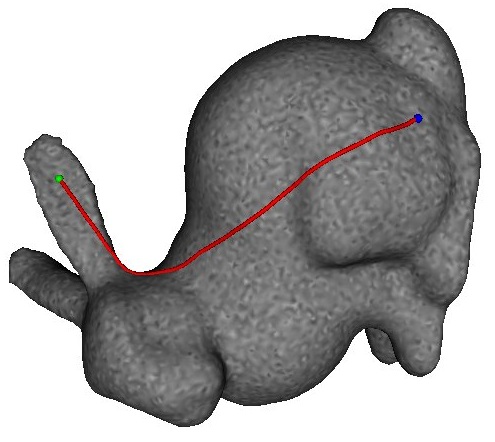}}

  \caption{\emph{Deblurring:} The first row shows the original signal and the noisy signal. The second row depicts the results with total variation regularization and the residual images, respectively. The last rows results are obtained with by regularizing with the squared gradient.}
  \label{fig:resultsDeblurring}
\end{figure}

Again, we compared TV and quadratic regularization. For TV regularization, $\kappa < \frac{2}{1+\alpha 8/ \epsilon}$ has to be chosen sufficiently
small.
Figure \ref{fig:resultsDeblurring} shows results for deblurring with TV minimization and quadratic Tikhonov regularization.

\begin{table}
\centering
\begin{tabular}{cc c c c}
\hline
 &Type          & SNR Original/Noise& Original/Result \\
\hline\hline

\multirow{2}{*}{Sobolev semi-norm} & & \\
 & Denoising  & 22.12	 &  23.87  \\
 & Deblurring & 16.89	 &   17.23  \\
 \hline
\multirow{2}{*}{TV semi-norm} & &  \\
 & Denoising  & 22.12	 &  26.35  \\
 & Deblurring & 16.89	 &   21.18  \\\hline
\hline
\end{tabular}
\caption{This table evaluates TV and quadratic variational regularization for denoising and deblurring in terms of the signal to noise ratio(c.f. \ref{eq:SNReq}).}
\label{tab:SNR}
\end{table}
In Table \ref{tab:SNR} we summarized the results on the denoising and delurring problem for a fixed $\alpha$. In order to compare the performance of different
choices for $\Rr$ we use the signal-to-noise ratio (SNR) measured in $dB$ . The SNR is defined as
\begin{equation}
\label{eq:SNReq}
\text{ SNR} = 20 \log_{10} \left( \frac{\norm{u^\dagger}_{\text{\upshape L}^2(M)}}{\norm{\uad-u^\dagger}_{\text{\upshape L}^2(M)}} \right).
\end{equation}
The better performance of the Total Variation regularization stems from the fact that discontinuities along edges are preserved while
the Sobolev semi-norm introduces severe blurring of the edges.

\subsubsection*{Variational regularization for inversion of the spherical Funk-Radon transform }
In a recent work Louis et al \cite{LouRipSpiSpo11} discuss a problem of density estimation, which requires the
inversion of the Funk-Radon transform on the $2$-sphere.
In general, for arbitrary space dimension, the Funk-Radon transform maps a function defined on the $2$-sphere to its
means over the great circles. That is,
\begin{equation}
\label{eq:FunkTransform}
F u(\xtt) = \frac{1}{\omega_1}  \int_{\mathbb S^{2} \cap \xtt^\perp} u(\ytt) d\nu(\ytt) \text{ for every } \xtt \in \mathbb S^{2}\;.
\end{equation}
Helgason \cite{Hel99} provids a closed form for the inverse of the Funk-Radon transform.
In \cite{LouRipSpiSpo11} an approximate inverse for the efficient numerical inversion
of the Radon-Funk transform on the $2$-sphere has been proposed.

Here we investigate quadratic Tikhonov regularization with Sobolev semi-norm regularization term on the $2$-sphere. The method
consists in minimization of the functional
\begin{equation}
\TTT(u) =  \norm{Fu - \yd}_{\text{\upshape L}^2 (\mathbb S^{2})}^2 + \alpha \norm{\nabla_M u}_{\text{\upshape L}^2 (\mathbb S^{2})}^2\;.
\label{eq:FunkRadon}
\end{equation}
The proposed numerical minimization algorithm requires real valued spherical harmonics : The functions $Y_{l}^{m}(\theta,\phi)$, where $l$
denotes the degree and $m$ the order, form an orthonormal basis on $\mathbb S^2$:
\begin{equation*}
\int_{\mathbb S^{2}} Y^m_\ell (\xtt) Y^{m'}_{\ell'}(\xtt)\,d\xtt =
\int_{\theta=0}^\pi \int_{\phi=0}^{2 \pi}  Y^m_\ell (\theta,\phi) Y^{m'}_{\ell'}(\theta,\phi)  d(\theta,\phi) = \delta_{\ell \ell'}\delta_{mm'}.
\end{equation*}
In the following we define a single index $j:=j(l,m):=(l+1)l+m $ for $l = 0,1,2,\dots,L$  and $m = 0,\ldots,l$
and identify the coordinates $\xtt$ on the sphere with polar coordinates $(\theta,\phi)$.
\begin{figure}
  \centering
 \subfloat[Approximation of $u^\dagger$ with spherical harmonics basis of maximal degree $26$.]{\includegraphics[scale=0.15]{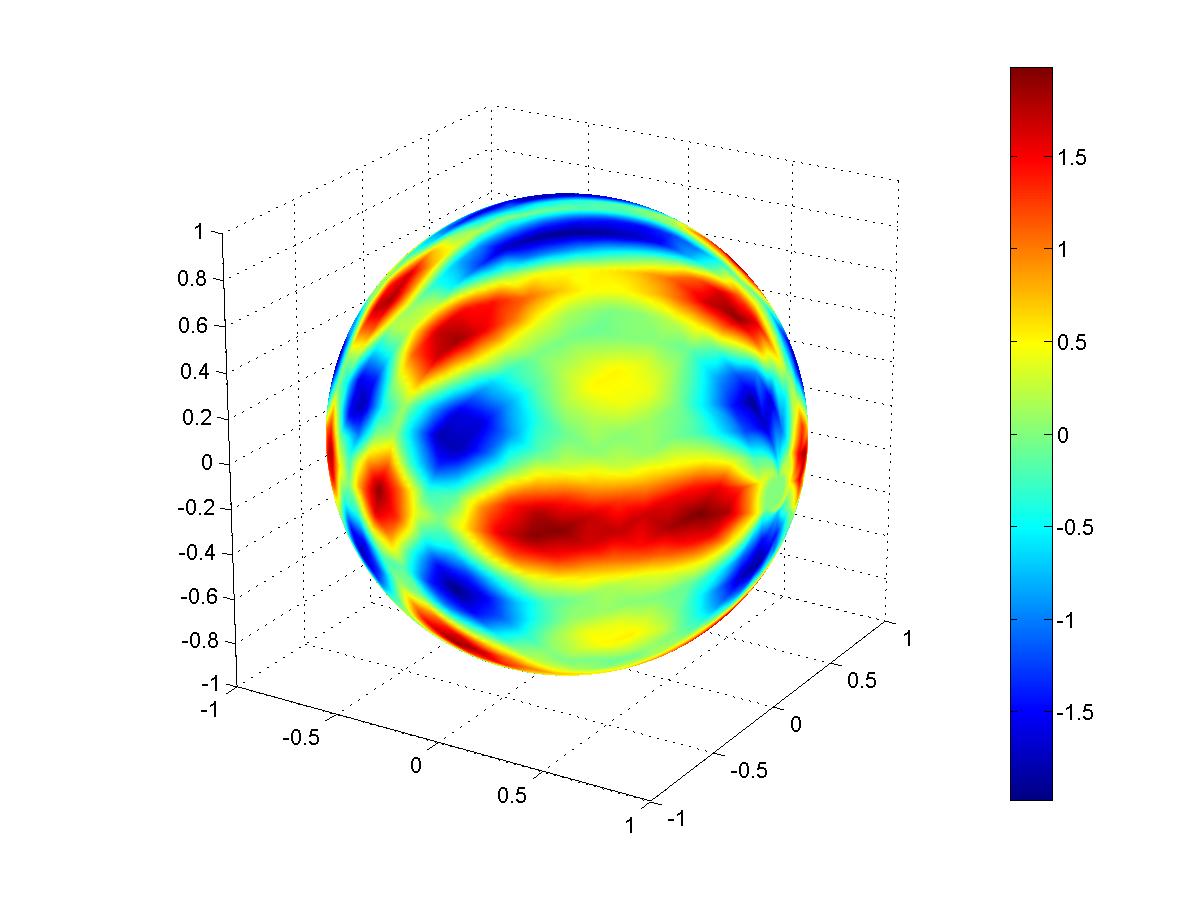}}
  \subfloat[Approximation error.]{\includegraphics[scale=0.15]{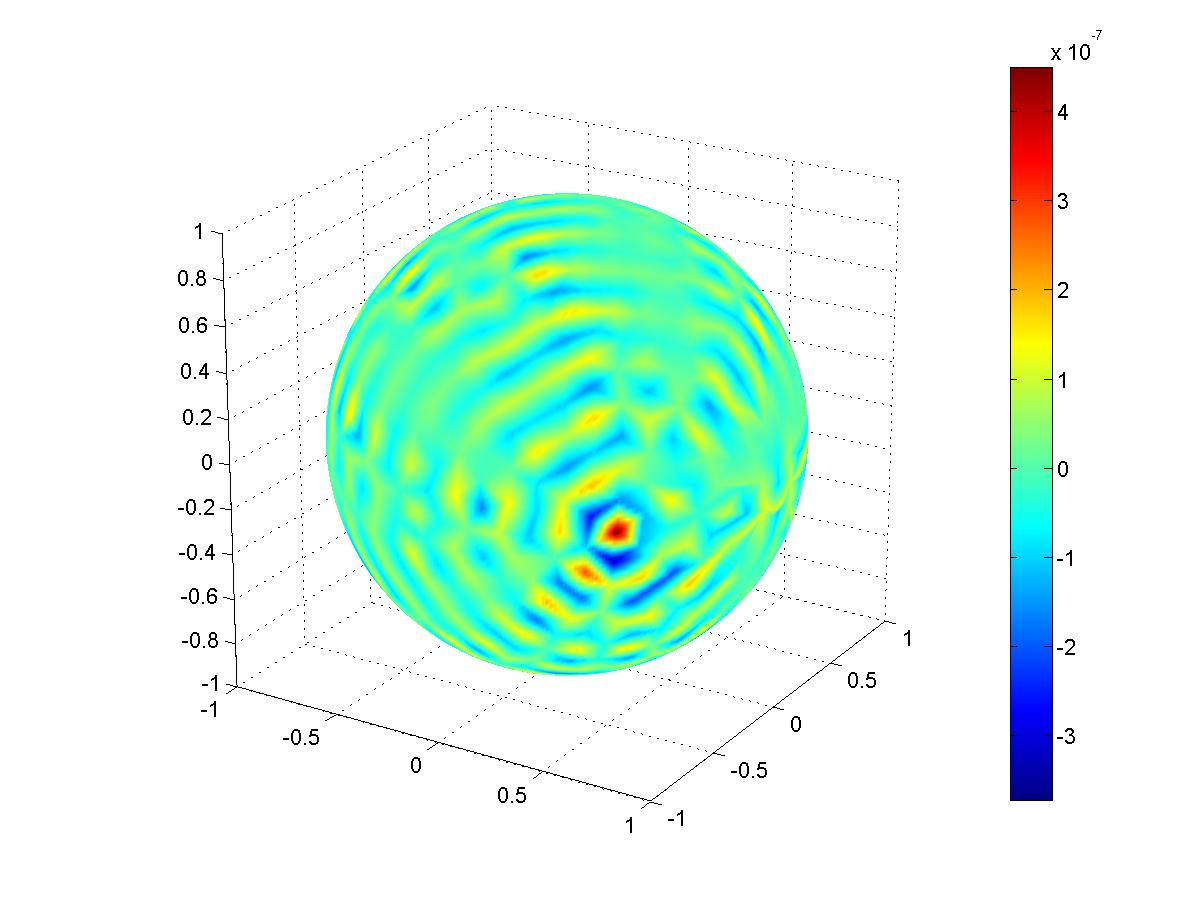}}

 \subfloat[$y=$ Funk-Radon transform of $u^\dagger$]{\includegraphics[scale=0.15]{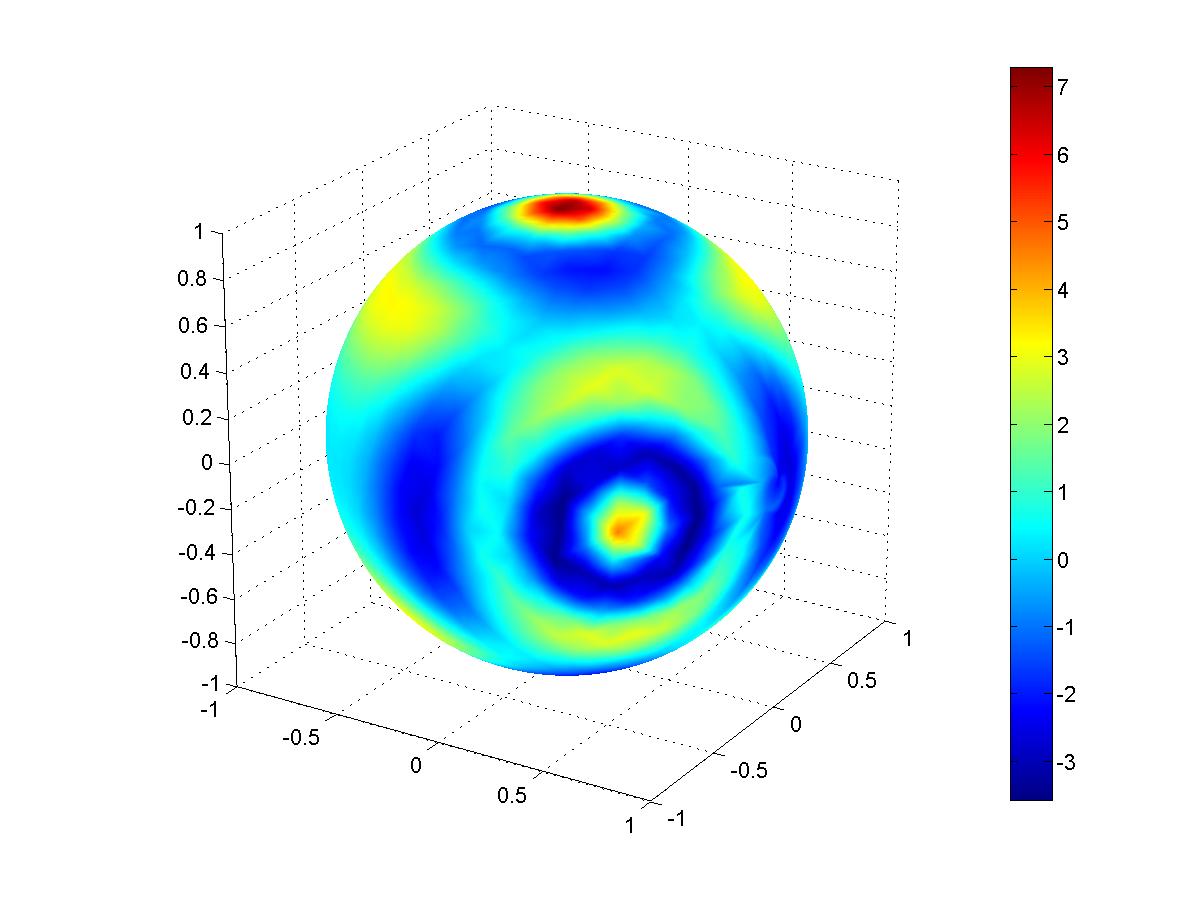}}
 \subfloat[$\yd =y$ + additive Gau\ss{}ian noise]{\includegraphics[scale=0.15]{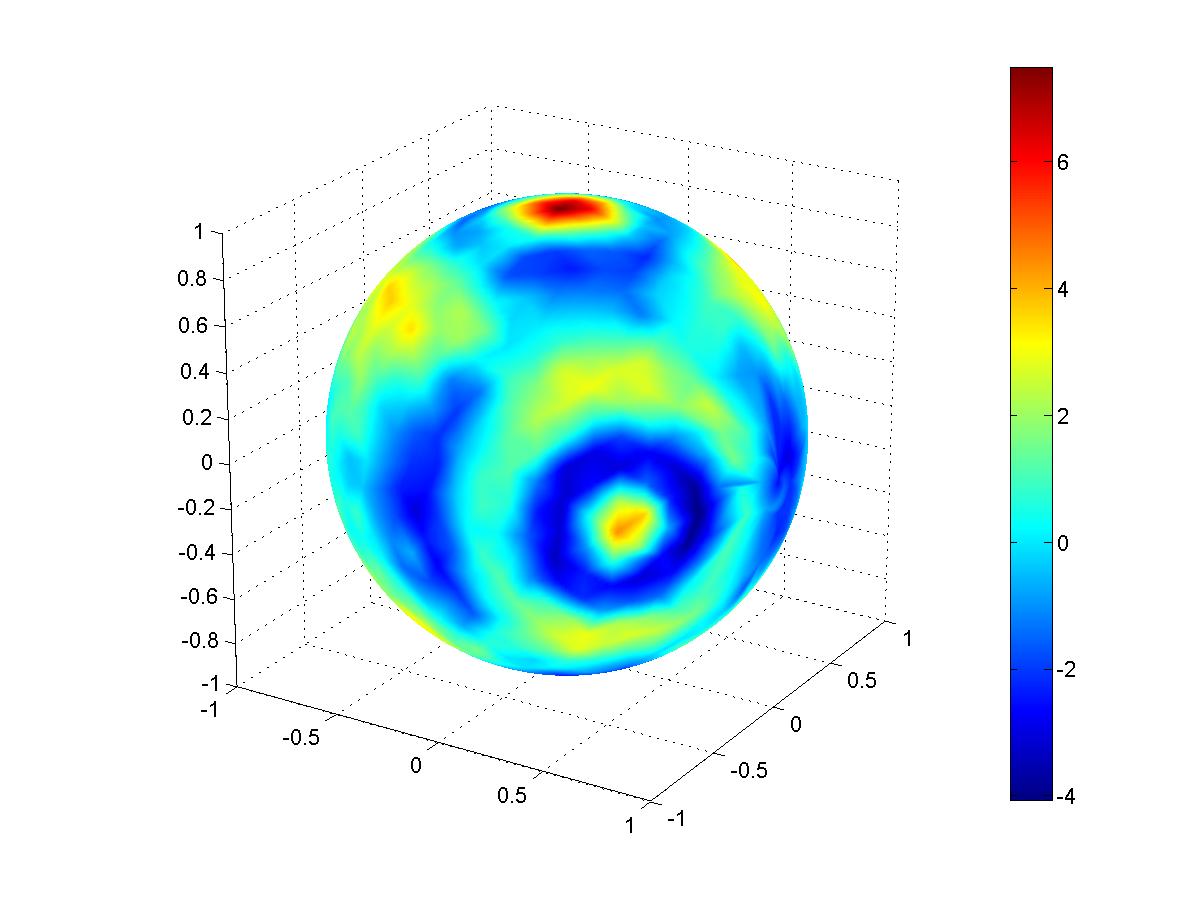}}
   \caption{\emph{Testdata and their Funk-Radon Transform} \label{fig:FunkTransform}}
\end{figure}
For the numerical minimization we use approximations of $u \in \text{\upshape L}^2(M)$ with real spherical harmonics of maximal
degree $L$. This is
\begin{equation}
u(\theta,\phi) \approx \sum_{j=1}^R  c_j Y_j(\theta,\phi).
\end{equation}
Therefore from discrete sample values $\mathbf u = (u_k)=(u(\xtt_k))\,,$ $k=1,\ldots,N \times R$ the spherical harmonics
expansions can be computed from the following matrix equation:
\[
\mathbf u = \mathbf B \mathbf c\,,
\]
where
\[
\mathbf B =	\left(
	\begin{array}{ccc}
		Y_1(\theta_1,\phi_1)& \dots & 		Y_R(\theta_1,\phi_1)\\
\vdots&\vdots &\vdots\\
		Y_1(\theta_N,\phi_N)& \dots & 		Y_R(\theta_N,\phi_N)\\
	\end{array} %
	\right) %
\]
is the matrix of spherical harmonics basis functions.
The coefficients $\mathbf c$ are the coefficients of the best approximating solution in $\text{\upshape L}^2(M)$ and are given by
\[
\mathbf c = (\mathbf B^T \mathbf B)^{-1}\mathbf B^T \mathbf u\;.
\]
Figure \ref{fig:FunkTransform} a) shows the best approximation of a function $u$ with a spherical harmonics polynomial of degree $26$.

Using the Funk-Henke Theorem it has been shown in \cite{DesAngFitDer06} that the Funk-Radon transform of a function
$u$, given in a spherical harmonics basis, takes the simple form
\begin{equation}
F u(\xtt) \approx   2\pi \sum_{j=1}^R P_{l(j)}(0) c_j Y_j(\xtt)\;.
\end{equation}
The Legendre polynomial of degree $l(j)$ evaluated at $0$ is
\begin{equation*}
P_{l(j)}(0)=
\begin{cases}
0 & \text{if $l(j)$ is odd,}\\
(-1)^{l(j)/2} \frac{1 \cdot 3 \cdot  \dots (l(j)-1)}{2 \cdot 4 \cdot 6 \dots l(j)} &\text{if $l(j)$ is even.}
\end{cases}
\end{equation*}
Therefore, the discrete Funk-Radon transform can be written as
\[
 \mathbf F = \mathbf B \hat{\mathbf F}\,, \text{ with } \hat{\mathbf F}_{jj} = 2\pi P_{l(j)}(0)
\]
a diagonal matrix.
An example of the evaluation of the Funk-Radon transform can be seen in Figure \ref{fig:FunkTransform} c).
In Figure \ref{fig:FunkTransform} d), the transformed signal is perturbed by some additive Gau\ss{}ian noise with variance
$\sigma= 0.05$.  In order to reconstruct the signal from its Funk-Radon transform we minimize the Tikhonov functional from Equation
(\ref{eq:FunkRadon}).

As in the previous examples the subgradient of the $\text{\upshape L}^2(M)$-norm of the gradient the Laplace-Beltrami operator, which in
spherical harmonics basis expansion is given by
\[
- \Delta_{\mathbb S^2} Y_j =  l(j)(l(j)+1) Y_j
\]
and in matrix notation
\[
\mathbf L^{jj} = l(j)(l(j)+1).
\]
Reconstruction of the inverse Funk-Radon transform with Tikhonov regularization requires solving the linear system
\begin{equation}
(\mathbf F^{\text{T}} \mathbf B^{\text{T}} \mathbf B \mathbf F + \alpha \mathbf L)\mathbf {\hat c} =
\mathbf F^{\text{T}} \mathbf B^{\text{T}} \mathbf \yd\;.
\label{invRadon}
\end{equation}
This equation can be solved again with a Landweber iteration. This equation has been solved again via a Landweber iteration:
\begin{equation}
\mathbf{\hat c}^{(n+1)} = \mathbf{\hat c}^{(n)} - \kappa \left(
\mathbf F^{\text{T}} \mathbf B^{\text{T}} (\mathbf B \mathbf F {\hat c^{(n)}} - \mathbf \yd)
+ \alpha \mathbf L\mathbf {\hat c}^{(n)} \right) \;.
\label{invRadon_Landweber}
\end{equation}

\begin{figure}
  \centering
	 \subfloat[Reconstruction with low regularization]{\includegraphics[scale=0.15]{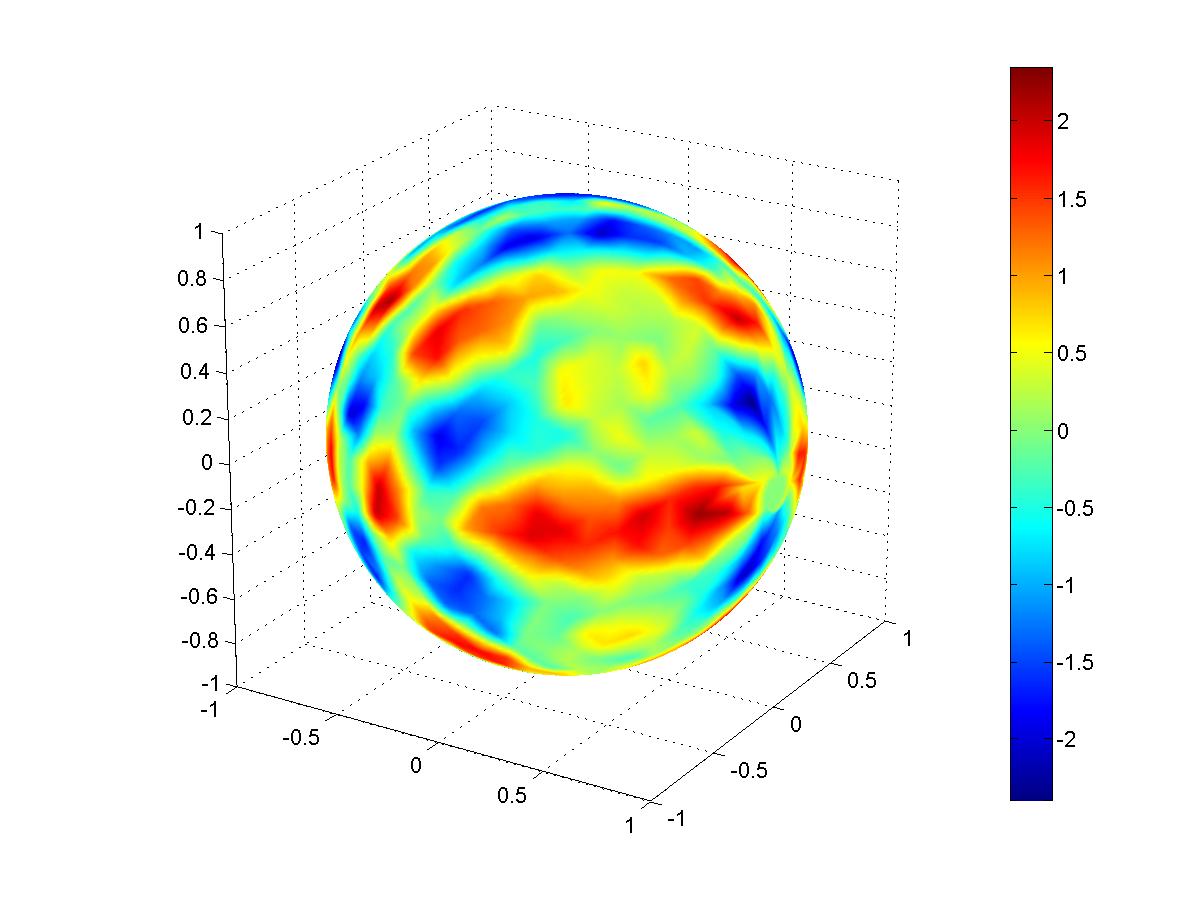}}
	 \subfloat[Reconstruction with strong regularization]{\includegraphics[scale=0.15]{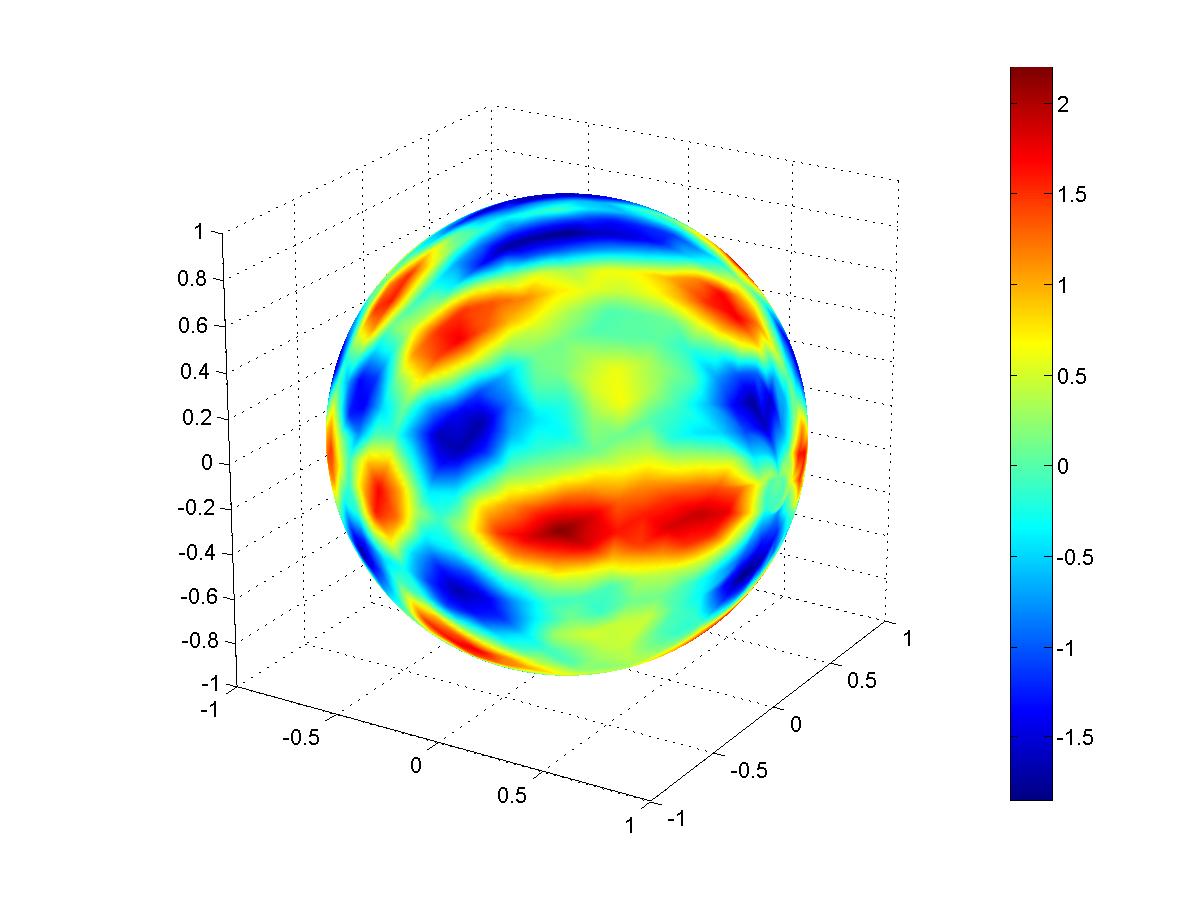}}

	 \subfloat[Reconstruction error with low regularization]{\includegraphics[scale=0.15]{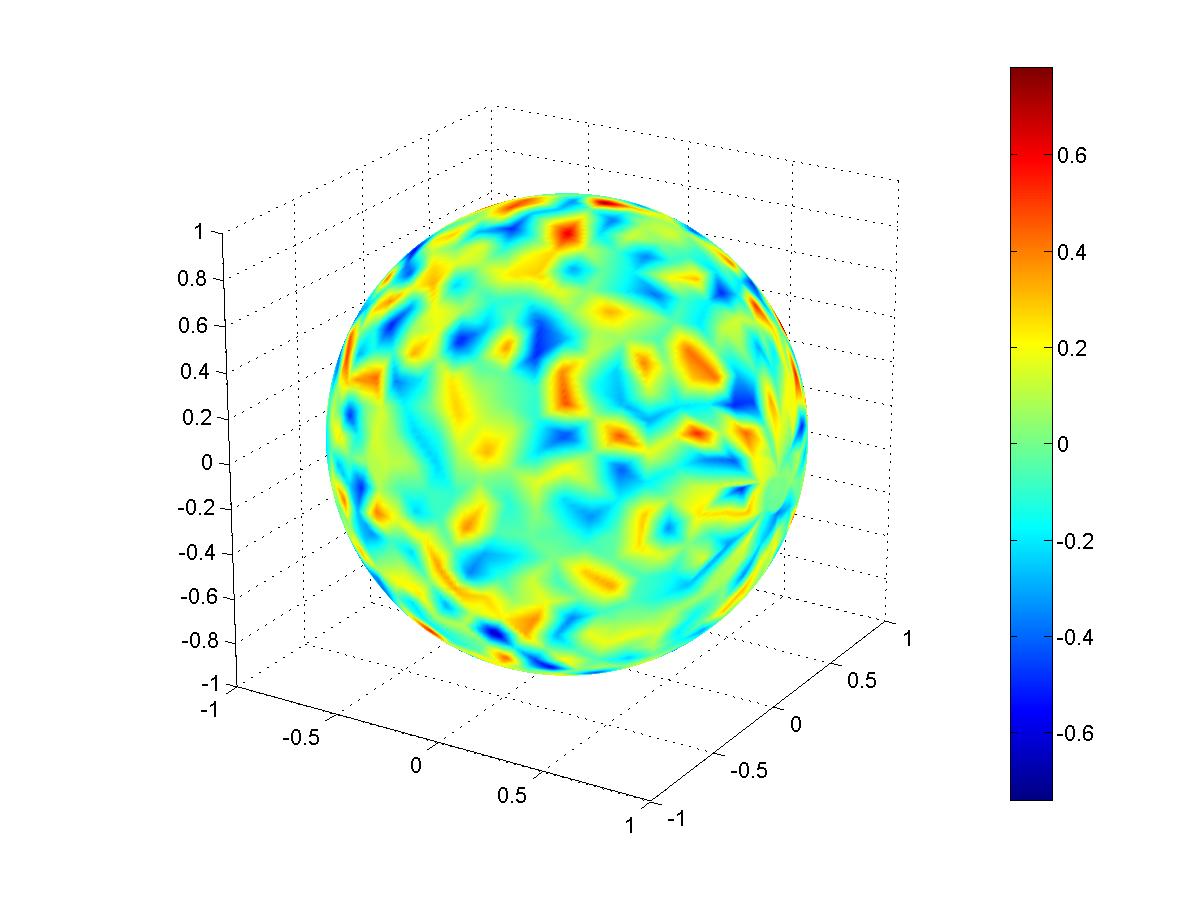}}
	 \subfloat[Reconstruction error with strong regularization]{\includegraphics[scale=0.15]{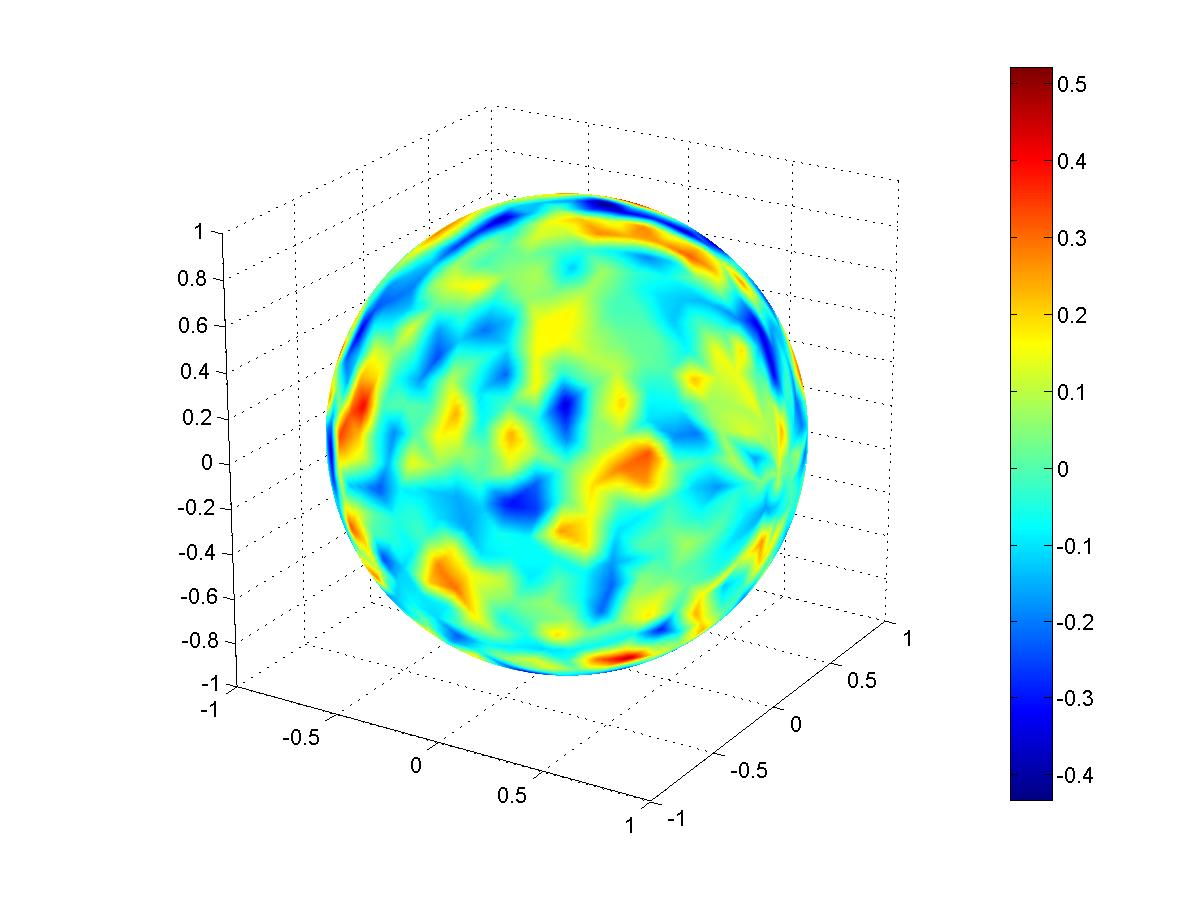}}
\caption{}
	\label{fig:InvFunkTransform}
\end{figure}

Since the Funk transform annihilates odd functions (see \cite{Hel99}), we take an even function to test our inversion algorithm.
As in \cite{LouRipSpiSpo11} we use the function
\[
u(\xtt) = \cos(3 \pi (\ztt-\ytt))+\cos(3 \pi \xtt);
\]
and evaluate the function at  $900$ point on the sphere as they are provides in \cite{Wom01}.
In the numerical experiments we used spherical harmonics of degree $26$. The reconstructions are depicted in Figure \ref{fig:InvFunkTransform}.
Note that if $\alpha =0$ we cannot solve for $\mathbf {\hat C} $ because the matrix $\mathbf F$ does not have full rank.
Only by regularizing the inversion with the Laplace-Beltrami operator allows the reconstruction of $u(\xtt)$. In Figure \ref{fig:InvFunkTransform} we can observe the smoothing effect of the Laplace-Beltrami operator on the solution. In the left column of Figure \ref{fig:InvFunkTransform} a low value $\alpha = 0.3$   and in the right column the result with  a high value $\alpha = 1.2$. We observe a much smoother reconstruction for higher values of $\alpha$.

\section{Conclusion}
\label{sec:conclusion}
In this paper, we have studied the problem of variational regularization of inverse and ill--posed problems for functions on closed Riemannian manifolds.
The analysis (stability, convergence, and rates) follows from standard results on convex regularization and are reviewed. The main contribution of this paper
concerns the numerical analysis of such regularization methods and the numerical implementation. Moreover, three inverse problems appearing in non-destructive evaluation and Computer Vision are discussed.

\section*{Acknowledgment}
The authors acknowledge the support by the Austrian Science Fund (FWF)
within the national research networks Industrial Geometry, project
9203-N12 and Photoacoustic Imaging, project P10505.


\section{Background and definitions}
\label{sec:def}
All along this paper we use the notation:
\begin{itemize}
\item $\Mm$ denotes a closed $n$-dimensional manifold in $\RR^{n+1}$.
\item $M = (\mathcal V,\mathcal T)$ denotes a polygonal approximation of $\Mm$. In this paper, typically, it is a polyhedron.
\item $\xtt$ denotes coordinates on the manifold. $x$ denotes coordinates in the Euclidean space.
\item $\partial$ denotes the subdifferential, $\nabla$ is the gradient in the Euclidean setting, $\nabla_M$ and $\nabla_\Mm$ are the covariant
derivatives on the Riemannian manifolds $M$, $\Mm$, respectively.
\item If not specified otherwise $\abs{.}$ denotes an Euclidean distance.
\end{itemize}

In the following we review elementary facts from Riemannian geometry and nonlinear analysis on manifolds.
\begin{enumerate}
\item The metric tensor $ g_{ij}(x)$ expressed in a coordinate chart $(\Om,\phi)$ is
\[
 G = g_{ij}(x) :=\inner{\frac{\partial \phi(x)}{\partial x_i}}{\frac{\partial \phi(x)}{\partial x_j}}_{i,j=1,\dots,n}
\]
	The inverse metric tensor is $g^{ij}=G^{-1}$.
\item Given a smooth, closed Riemannian $n$-manifold $(\Mm,g)$, there is an associated positive Radon measure on $\Mm$, the \emph{Riemannian measure}, which is defined as follows: Given an integrable function $u : \Mm \mapsto \RR$,
    an atlas $(\Omega_i,\phi_i)_{i \in I}$ of $\Mm$ and a partition of unity  $(\Omega_j, \phi_j, \alpha_j)_{j \in J}$,
\begin{equation}
\label{eq:manifoldIntegral}
    \int_{\Mm} u d \nu(g) = \sum_{j\in J} \int_{\phi_j( \Omega_j)} (\alpha_j  u\circ \phi_j \sqrt{\abs{g}} )  dx\,,
\end{equation}
    where $\sqrt{\abs{g}} = |\text{det}(g_{ij})|$ is the volume of the metric tensor and $d x$ is the Lebesgue volume element on $\RR^n$ and therefore $d \nu(g) =\sqrt{\abs{g}}   dx $
\item The gradient can be expressed with the chart $(\Om,\phi)$
\begin{equation}
\label{eq:nabla}
(\nabla_\Mm u)_i =  g_{ij}  \frac{\partial( u \circ \phi^{-1})}{\partial x_j} \;.
\end{equation}
\item  The adjoint operator of the gradient is the divergence, which satisfies for a given  vector field $X$ on $\Mm$:
\begin{align}
\label{eq:divM}
\int_\Mm u \div_\Mm X d \nu(g) =- \int_{\Mm} \nabla_\Mm  u \cdot X d \nu(g)\;.
\end{align}
In a chart $(\Omega,\phi)$ the divergence is obtained by using Equation \ref{eq:nabla} in Equation \ref{eq:divM}:
\begin{equation}
\label{eq:divergence}
\div_\Mm X = \frac{1}{\sqrt {\abs{g}}}  \sum_{i=1}^n  \frac{\partial}{\partial x_i}(\sqrt {\abs{g}}X_i).
\end{equation}
\item Given $(\Mm,g)$ and $\gamma :[ a,b ]\mapsto \Mm$ a curve on $\Mm$, then $L(\gamma)$ is the length of the curve on $\Mm$ with respect to $g$. For $p,q $ on $\Mm$ with $\gamma(a)=p$ $\gamma(a)=p$ and $\gamma(b)=q$, the distance associated with $g$ between two points $p$ and $q$ is
    \[
        d_g(p,q) = \inf_{\gamma \in C_{pq}}L(\gamma).
    \]
$C_{pq}$ is the space of continuous curves connecting $p$ and $q$. The distance $d_g$ defines a metric on the manifold. In this paper, we assume that the metric space $(\Mm,d_g)$ is always complete.
\item In a \emph{closed Riemannian manifold} (without boundary) the Hopf-Rinow theorem implies that for every pair of points on the manifold there exists a unique geodesic \cite{Cha84}.
\end{enumerate}

Given a smooth and closed $n$-dimensional Riemannian manifold $(\Mm,g)$, we define
(see \cite{Heb96})
\begin{equation*}
\begin{aligned}
& \Cc^{p}_{k}(\Mm ) := \\
&\set{u \in  \Cc^{\infty}(\Mm) : \norm{u}_{k,p} := \norm{u}_p +\sum_{j=1}^{k}
\abs{\int_{\Mm} \abs{\nabla_\Mm^j u}^p d \nu(g)}^{\frac{1}{p}}  < \infty}\,,
\end{aligned}
\end{equation*}
where $\norm{.}_{p}$ is the $\text{\upshape L}^p(\Mm)$-norm with respect to the Riemannian measure $d\nu(g)$.
The space $W^{p,k}(\Mm)$ (see e.g. \cite{Heb96}) is defined as the completion of the space
$\Cc^{p}_{k}(\Mm )$ with respect to the norm $\norm{.}_{k,p}$. In particular for $k = 1$ and $p > 1$, we have
\[
\norm{u}_{W^{p,1}(\Mm)} = \norm{u}_{p} + \norm{\nabla_\Mm u}_{p}\;.
\]
Now, we recall the definition of the space of functions of bounded variation on manifolds.
\begin{defi}
We define $BV(\Mm)$ as the space of functions with bounded variation and is the set of functions $u \in \text{\upshape L}^1(\Mm)$ such that
$\abs{D_\Mm u}(\Mm) < + \infty$. The space is endowed with the
norm $\norm{u}_{BV(\Mm)} = \norm{u}_{\text{\upshape L}^1(\Mm)} + \abs{D_\Mm u}(\Mm)$\,,
where $\abs{D_\Mm u}(\Mm)$ denotes the variation of $u$, which is defined by (\ref{eq:tv})
\end{defi}
The space $BV(\Mm)$ is a Banach space endowed with the norm $\norm{.}_{BV(\Mm)}$. It can be understood as the natural (weak) closure of
$W^{1,1}(\Mm)$. Due to the theorem of Meyer and Serrin \cite{MeySer64} it is possible to approximate Sobolev functions defined on subset
of the Euclidean space by smooth functions.
For the sake of completeness we provide a proof to show the essential difference in the manifold setting.
\begin{theo}[Approximation of $BV$-Functions]
\label{meyerSerr}
Let $\Mm$ be a smooth, closed Riemannian manifold and $u \in BV(\Mm)$. Then there exists a sequence $(u_n)$ of functions in $\Cinf_c(\Mm)$ such that
\begin{itemize}
\item $u_n \rightarrow u$ in $\text{\upshape L}^1(\Mm)$,
\item $\int_{\Mm} \abs{\nabla_\Mm u_n(\xtt)} d\xtt \rightarrow \abs{D_\Mm u}(\Mm)$.
\end{itemize}
\end{theo}
\begin{prof}
The proof is closely related to \cite{AmbFusPal00}, where weighted $BV$ spaces have been considered and thus omitted.
\end{prof}
The second important property of $BV$ functions used in this paper is covered by the following embedding theorem:
\begin{theo}[Compactness Theorem]
\label{th:compactness}
Let $\Mm $ be a closed manifold, and let $(u_n)_n$ be a sequence of functions in $BV(\Mm)$ such that $\sup_n \abs{D_\Mm u_n}(\Mm) \leq + \infty$. Then there exists a subsequence of $u \in BV(\Mm)$ converging strongly in $\text{\upshape L}^1(\Mm)$.
\end{theo}
Follows from combining the analogous result for functions $W^{1,1}(\Mm)$, which is stated in Hebey\cite{Heb96}, and Theorem \ref{meyerSerr}.

\begin{theo}[Embedding Theorem]
\label{th:embedding}
For every function $u \in BV(\Mm)$
\begin{equation}
\label{eq:sobolevEQ}
\norm{u}_{\text{\upshape L}^{\frac{n}{n-1}}(\Mm)} \leq C(n) \abs{D_\Mm u}(\Mm).
\end{equation}
\end{theo}
The proof is analogous to the Euclidean setting and thus omitted.


\def\cprime{$'$} \providecommand{\noopsort}[1]{}

\end{document}